\documentclass[10pt]{amsart}
\usepackage{amsfonts}
\usepackage{ifthen}
\usepackage{amsthm}
\usepackage{amsmath}
\usepackage{graphicx}
\usepackage{amscd,amssymb,amsthm}
\usepackage{graphicx}
\usepackage{epstopdf}
\usepackage{hyperref}

\newcounter{minutes}
\divide\time by 60
\newcounter{hours}
\multiply\time by 60 \addtocounter{minutes}{-\time}

\newtheorem{theorem}{Theorem}

\newcommand{\real}{\operatorname{Re}}


\begin{document}

\title{Bounds for radii of starlikeness and convexity of some special functions}

\author[\.{I}. Akta\c{s}]{\.{I}brah\.{i}m Akta\c{s}}
\address{Department of Mathematical Engineering, Faculty of Engineering and Natural Sciences, G\"{u}m\"{u}\c{s}hane University, G\"{u}m\"{u}\c{s}hane, Turkey}
\email{aktasibrahim38@gmail.com}
\author[\'A. Baricz]{\'Arp\'ad Baricz$^{\bigstar}$}
\address{Department of Economics, Babe\c{s}-Bolyai University, Cluj-Napoca, Romania}
\address{Institute of Applied Mathematics, \'Obuda University, Budapest, Hungary}
\email{bariczocsi@yahoo.com}
\author[H. Orhan]{Hal\.{i}t Orhan}
\address{Department of Mathematics, Faculty of Science, Atat\"{u}rk University, Erzurum, Turkey}
\email{orhanhalit607@gmail.com}

\thanks{$^{\bigstar}$The research of \'A. Baricz was supported by a research grant of the Babe\c{s}-Bolyai University for young researchers with project number GTC-31777.}
\maketitle

\begin{abstract}
In this paper we consider some normalized Bessel, Struve and Lommel functions of the first kind, and by using the Euler-Rayleigh inequalities for the first positive zeros of combination of special functions we obtain tight lower and upper bounds for the radii of starlikeness of these functions. By considering two different normalization of Bessel and Struve functions we give some inequalities for the radii of convexity of the same functions. On the other hand, we show that the radii of univalence of some normalized Struve and Lommel functions are exactly the radii of starlikeness of the same functions. In addition, by using some ideas of Ismail and Muldoon we present some new lower and upper bounds for the zeros of derivatives of some normalized Struve and Lommel functions. The Laguerre-P\'olya class of real entire functions plays an important role in our study.
\end{abstract}

\section{Introduction}
It is known that special functions, like Bessel, Struve and Lommel functions of the first kind and regular Coulomb wave function have some beautiful geometric properties. Recently, the geometric properties of the above special functions were investigated motivated by some earlier results. In the sixties Brown, Kreyszig and Todd, Wilf (see \cite{brown,brown2,brown3,todd,wilf}) considered the univalence and starlikeness of Bessel functions of the first kind, while in the recent years the radii of univalence, starlikeness and convexity for the normalized forms of Bessel, Struve and Lommel functions of the first kind were obtained, see the papers \cite{aktas,mathematica,publ,bsk,bos,samy,basz,basz2,BY,szasz,szasz2} and the references therein. In these papers it was shown that the radii of univalence, starlikeness and convexity are actually solutions of some transcendental equations. On the other hand, it was shown that the obtained radii satisfy some interesting inequalities. In addition, it was proved that the radii of univalence of some normalized Bessel and Struve functions correspond to the radii of starlikeness of the same functions. In the above works the authors used intensively some properties of the positive zeros of Bessel, Struve and Lommel functions of the first kind, under some conditions. Also, they utilized the Laguerre-P\'{o}lya class $\mathcal{LP}$ of real entire functions. Motivated by the above developments in this topic, in this paper our aim is to give some new results for the radii of univalence, starlikeness and convexity of the normalized Bessel, Struve and Lommel functions of the first kind. This paper is a direct continuation of the paper \cite{aktas} and it is organized as follows: section 1 contains some basic concepts, in section 2 we focus on linear combination of Struve function and its derivative and the derivative of Lommel function. Here we give some lower and upper bounds for the smallest positive zeros of these functions. To prove our results we use some ideas from \cite{ismail}. We also consider two normalized forms of Struve and Lommel functions, respectively. For these functions, we show that the radii of univalence and starlikeness coincide. At the end of this section we obtain some new lower and upper bounds concerning the radii of convexity of four different normalized forms of Bessel and Struve functions of the first kind. In section 3 we present the proofs of the main results. The bounds deduced for the radii of convexity are in fact particular cases of some Euler-Rayleigh inequalities and it is possible to show that the lower bounds increase and the upper bounds decrease to the corresponding radii of convexity, and thus the inequalities presented in this paper can be improved by using higher order Euler-Rayleigh inequalities. We restricted ourselves to the third Euler-Rayleigh inequalities since these
are already complicated.

Now, we would like to present some basic concepts regarding geometric function theory. Let $\mathbb{D}_r=\{z\in\mathbb{C}:\left|z\right|<r\}$ be the open disk, where $r>0$. Also, let $f:\mathbb{D}_r\rightarrow\mathbb{C}$ be the function, defined by
\begin{equation}
f(z)=z+\sum_{n\geq2}a_{n}z^{n}.
\end{equation}

The function $f$, defined by $(1.1)$, is called starlike in the disk $\mathbb{D}_r$ if $f$ is univalent in $\mathbb{D}_r$, and $f(\mathbb{D}_r)$ is a starlike domain in $\mathbb{C}$ with respect to origin. Analytically, the function $f$ is starlike in $\mathbb{D}_r$ if and only if $$\real\left(\frac{zf^{\prime}(z)}{f(z)}\right)>0 \text{ for all } z\in\mathbb{D}_r.$$ The real number $$r^{\ast }(f)=\sup \left\{ r>0 \left|\real%
\left(\frac{zf^{\prime }(z)}{f(z)}\right) >0 \;\text{for all }z\in
\mathbb{D}_r\right.\right\}$$ is called the radius of starlikeness of the function $f$.

The function $f,$ defined by $(1.1)$, is convex in the disk $\mathbb{D}_r$ if $f$ is univalent in $\mathbb{D}_r$, and $f(\mathbb{D}_r)$ is a
convex domain in $\mathbb{C}.$ Analytically, the function $f$ is convex in $\mathbb{D}_r$ if and only if $$\real\left(1+\frac{zf^{\prime \prime}(z)}{f^{\prime}(z)}\right)>0 \text{ for all } z\in\mathbb{D}_r.$$ The radius of convexity of the function $f$ is defined by the real number $$r^{c}(f)=\sup \left\{r>0 \left|\real\left(1+\frac{zf^{\prime \prime}(z)}{f^{\prime}(z)}\right)>0 \;\text{ for all }
z\in\mathbb{D}_r\right.\right\}.$$
Finally, we recall that the radius of univalence of the analytic function $f$ in the form of $(1.1)$ is the largest radius $r$ such that $f$ maps $\mathbb{D}_r$ univalently into $f(\mathbb{D}_r).$

\section{Bounds for the zeros of some special functions}
In this paper we consider three classical special functions, the Bessel function of the first kind $J_\nu$, the Struve function of the first kind $\mathbf{H_\nu}$ and the Lommel function of the first kind $s_{\mu,\nu}$. It is known that the Bessel functions has the infinite series representation $$J_{\nu}(z)=\sum_{n\geq0}\frac{(-1)^n}{n!\Gamma(n+\nu+1)}
\left(\frac{z}{2}\right)^{2n+\nu},$$ where $z,\nu\in\mathbf{C}$ such that $\nu\neq-1,-2,{\dots}.$ Also, the Struve and Lommel functions can be represented as the infinite series $$\mathbf{H_\nu}(z)=\sum_{n\geq0}\frac{(-1)^n}{\Gamma\left(n+\frac{3}{2}\right)\Gamma\left(n+\nu+\frac{3}{2}\right)}\left(\frac{z}{2}\right)^{2n+\nu+1}, -\nu-\frac{3}{2}\notin\mathbb{N},$$ and $$s_{\mu,\nu}=\frac{(z)^{\mu+1}}{(\mu-\nu+1)(\mu+\nu+1)}\sum_{n\geq0}\frac{(-1)^n}{(\frac{\mu-\nu+3}{2})_n(\frac{\mu+\nu+3}{2})_n}\left(\frac{z}{2}\right)^{2n}, \frac{1}{2}(-\mu\pm\nu-3)\notin\mathbb{N},$$ where   $z,\mu,\nu\in\mathbf{C}$. In addition, we know that the Bessel function is a solution of the homogeneous Bessel differential equation $$zw^{\prime \prime}(z)+zw^{\prime}(z)+(z^2-\nu^2)w(z)=0,$$ while the Struve and Lommel functions are solutions of the inhomogeneous Bessel differential equations $$zw^{\prime \prime}(z)+zw^{\prime}(z)+(z^2-\nu^2)w(z)=\frac{4\left(\frac{z}{2}\right)^{\nu+1}}{\sqrt{\pi}\Gamma\left(\nu+\frac{1}{2}\right)}$$ and $$zw^{\prime \prime}(z)+zw^{\prime}(z)+(z^2-\nu^2)w(z)=z^{\mu+1},$$respectively. We refer to Watson's treatise \cite{Wat} for comprehensive information about these functions. On the other hand, the Laguerre-P\'{o}lya class $\mathcal{LP}$ of real entire functions plays an important role in our proofs. Recall that a real entire function $\Psi$ belongs to the  Laguerre-P\'{o}lya class $\mathcal{LP}$ if it can be represented in the form $$\Psi(x)=cx^{m}e^{-ax^2+bx}\prod_{n\geq1}\left(1+\frac{x}{x_n}\right)e^{-\frac{x}{x_n}},$$ with $c,b,x_n\in\mathbb{R}, a\geq0, m\in\mathbb{N}_0$ and $\sum1/{x_n}^2<\infty.$

We note that the class $\mathcal{LP}$ consists of entire functions which are uniform limits on the compact sets of the complex plane of polynomials with only real zeros. For more details on the class $\mathcal{LP}$ we refer to \cite[p. 703]{DC} and to the references therein.

\subsection{Zeros of linear combination of Struve function and its derivative}
In this subsection by considering the Struve function $\mathbf{H}_{\nu}$ and its derivative $\mathbf{H}_{\nu}^{\prime}$ we define the function $\mathcal{H}_{\nu}$ as follows:
$$\mathcal{H}_{\nu}(z) =\alpha \mathbf{H}_{\nu}+z\mathbf{H}_{\nu}^{\prime}(z).$$ The function $\mathcal{H}_{\nu}$ can be written as $$\mathcal{H}_{\nu}(z)=\sum\limits_{n\geq 0}\frac{(-1)^{n}(2n+{\nu}+\alpha +1) }{\Gamma (n+\frac{3}{2}) \Gamma({\nu}+n+\frac{3}{2})}\left(\frac{z}{2}\right)^{2n+{\nu}+1}.$$
Let $\alpha +\nu\neq -1.$ Here we focus on the following normalized form:
$$h_{\nu}(z)=(\alpha+{\nu}+1)^{-1}\Gamma\left(\frac{3}{2}\right)\Gamma\left({\nu}+\frac{3}{2}\right) z^{-\frac{{\nu}+1}{2}}2^{{\nu}+1}\mathcal{H}_{\nu}(\sqrt{z})=\sum\limits_{n\geq0}\frac{(-1))^{n}(2n+{\nu}+\alpha+1)}{2^{2n}({\nu}+\alpha+1)(\frac{3}{2})_{n}({\nu}+\frac{3}{2})_{n}}z^{n}$$ Our first main result is related to the function $\mathcal{H}_{\nu}.$
\begin{theorem}
	Let $\alpha +{\nu}> -1, \left|{\nu}\right|<\frac{1}{2}$ and let $\zeta_{{\nu},1}$ be the smallest positive zero of the function $\mathcal{H}_{\nu}.$ Then we have the lower bounds $$\zeta_{{\nu},1}^{2}>\frac{3(2{\nu}+3)(\alpha+{\nu}+1)}{\alpha+{\nu}+3},$$$$\zeta_{{\nu},1}^{2}>\frac{3(2{\nu}+3)(\alpha+{\nu}+1)\sqrt{5(2\nu+5)}}{\sqrt{\kappa_1}},$$$$\zeta_{{\nu},1}^{2}>\frac{3(2{\nu}+3)(\alpha+{\nu}+1)\sqrt[3]{35(2\nu+5)(2\nu+7)}}{\sqrt[3]{\kappa_2}}$$ and the upper bounds $$\zeta_{{\nu},1}^{2}<\frac{15(2\nu+3)(2\nu+5)(\alpha+\nu+1)(\alpha+\nu+3)}{\kappa_1},$$$$\zeta_{{\nu},1}^{2}<\frac{21(2\nu+3)(2\nu+7)(\alpha+\nu+1)\kappa_1}{\kappa_2},$$ where $\kappa _1=-2\alpha ^2\nu+7\alpha ^2-4\alpha\nu^2+2\alpha
	\nu+42\alpha -2\nu^3-5\nu^2+72\nu+135$ and $\kappa _2=-4\alpha
	^3\nu^2-96\alpha ^3\nu+145\alpha ^3-12\alpha ^2\nu^3-324\alpha
	^2\nu^2-429\alpha ^2\nu+1305\alpha ^2-12\alpha \nu^4-360\alpha
	\nu^3-1689\alpha \nu^2+1170\alpha\nu+6291\alpha
	-4\nu^5-132\nu^4-1115\nu^3+621\nu^2+12339\nu+14931.$
\end{theorem}
In particular, when $\alpha=0$, Theorem 1 reduces to the following.
\begin{theorem}
	Let $\left|{\nu}\right|<\frac{1}{2}$ and let $h_{\nu,1}^{\prime}$ be the smallest positive root of $\mathbf{H}_{\nu}^{\prime}$. Then we have the lower bounds $$(h_{{\nu},1}^{\prime})^2>\frac{3(2{\nu}+3)({\nu}+1)}{{\nu}+3},$$$$(h_{{\nu},1}^{\prime})^2>\frac{3(2{\nu}+3)({\nu}+1)\sqrt{5(2\nu+5)}}{\sqrt{-2\nu^3-5\nu^2+72\nu+135}},$$$$(h_{{\nu},1}^{\prime})^2>\frac{3(2{\nu}+3)({\nu}+1)\sqrt[3]{35(2\nu+5)(2\nu+7)}}{\sqrt[3]{-4\nu^5-132\nu^4-1115\nu^3+621\nu^2+12339\nu+14931}}$$ and the upper bounds $$(h_{{\nu},1}^{\prime})^2<\frac{15(2\nu+3)(2\nu+5)(\nu+1)(\nu+3)}{-2\nu^3-5\nu^2+72\nu+135},$$$$(h_{{\nu},1}^{\prime})^2<\frac{21(2\nu+3)(2\nu+7)(\nu+1)(-2\nu^3-5\nu^2+72\nu+135)}{-4\nu^5-132\nu^4-1115\nu^3+621\nu^2+12339\nu+14931}.$$
\end{theorem}
Here it is worth to mention that Theorem 2 reobtains and improves some results of \cite{BKP} regarding the first positive zeros of derivative of the Struve function. We mention that our approach is a little bit different than the approach in \cite{BKP}.

\subsection{Bounds for the zeros of derivative of Lommel functions}
We consider the function $$\mathcal{L}_\mu(z)=zs_{\mu-\frac{1}{2},\frac{1}{2}}^{\prime}(z)=\sum_{n\geq0}\frac{(-1)^n(2n+\mu+\frac{1}{2})}{4^n\mu(\mu+1)(\frac{\mu+2}{2})_{n}\frac{\mu+3}{2})_{n}}z^{2n+\mu+\frac{1}{2}},$$ where $s_{\mu-\frac{1}{2},\frac{1}{2}}^{\prime}(z)$ stands for the derivative of Lommel function. Let $\mu\in(-1,1), \mu\neq0$ and $\mu\neq-\frac{1}{2}.$ Now, we define the following normalized form of the function $\mathcal{L}_\mu$. Let $$l_{\mu}(z)=\frac{2\mu(\mu+1)}{(2\mu+1)}z^{-\frac{2\mu+1}{4}}\mathcal{L}_\mu(\sqrt{z}).$$ Clearly, the function $l_\mu$ can be written as $$l_{\mu}(z)=1+\sum_{n\geq1}{\frac{(-1)^n(2n+\mu+\frac{1}{2})}{2^{2n}(\mu+\frac{1}{2})(\frac{\mu+2}{2})_n(\frac{\mu+3}{2})_n}}z^n.$$
\begin{theorem}
	Let $\mu\in(-1,1), \mu\neq0, \mu\neq-\frac{1}{2}$ and let $\tau_{\mu,1}$ be the smallest positive zero of the function $\mathcal{L}_\mu$. Then we have the lower bounds
	$$(\tau_{\mu,1})^2>\frac{(\mu+2)(\mu+3)(2\mu+1)}{2\mu+5},$$
	$$(\tau_{\mu,1})^2>\frac{(\mu+2)(\mu+3)(2\mu+1)\sqrt{(\mu+4)(\mu+5)}}{\sqrt{-4\mu^4-24\mu^3+19\mu^2+295\mu+392}},$$
	 $$(\tau_{\mu,1})^2>\frac{(\mu+2)(\mu+3)(2\mu+1)\sqrt[3]{(\mu+4)(\mu+5)(\mu+6)(\mu+7)}}{\sqrt[3]{8\mu^7+44\mu^6-554\mu^5-4731\mu^4-7672\mu^3+23551\mu^2+85834\mu+72384}}$$ and the upper bounds
	$$(\tau_{\mu,1})^2<\frac{(\mu+2)(\mu+3)(\mu+4)(\mu+5)(2\mu+1)(2\mu+5)}{-4\mu^4-24\mu^3+19\mu^2+295\mu+392},$$
	 $$(\tau_{\mu,1})^2<\frac{(\mu+2)(\mu+3)(\mu+6)(\mu+7)(2\mu+1)(-4\mu^4-24\mu^3+19\mu^2+295\mu+392)}{8\mu^7+44\mu^6-554\mu^5-4731\mu^4-7672\mu^3+23551\mu^2+85834\mu+72384}.$$
\end{theorem}

\subsection{Radii of univalence (and starlikeness) of Struve functions}
Here our aim is to show that the radii of univalence of the Struve function $u_{\nu}$ correspond to the radii of starlikeness.
\begin{theorem}
	Let $\nu\in[-\frac{1}{2},\frac{1}{2}]$. The radius of univalence $r^*(u_\nu)$ of the normalized Struve function
	$$z\mapsto{u}_\nu(z)=\left(\sqrt{\pi}2^{\nu}\Gamma\left(\nu+\frac{3}{2}\right)\mathbf{H}_{\nu}(z)\right)^{\frac{1}{\nu+1}}$$
	corresponds to its radius of starlikeness and it is the smallest positive root ${{h}_{\nu,1}^{\prime}}$ of ${\mathbf{H}_{\nu}^{\prime}}.$
\end{theorem}

\subsection{Radii of univalence (and starlikeness) of Lommel functions}
In this subsection our aim is to show that the radii of univalence of the Lommel function $f_{\mu}$ correspond to the radii of starlikeness.
\begin{theorem}
	Let $\mu\in(-\frac{1}{2},1), \mu\neq0.$ The radius of univalence $r^*(f_\mu)$ of the normalized Lommel function $$z\mapsto{f}_\mu(z)={f}_{\mu-\frac{1}{2},\frac{1}{2}}(z)=\left(\mu({\mu+1}){s}_{\mu-\frac{1}{2},\frac{1}{2}}(z)\right)^{\frac{1}{\mu+\frac{1}{2}}}$$
	corresponds to its radius of starlikeness and it is the smallest positive root of ${s}_{\mu-\frac{1}{2},\frac{1}{2}}^{\prime}.$
\end{theorem}

\subsection{Radii of convexity of Bessel functions}
In this subsection we consider two different normalized forms of the Bessel functions of the first kind. Here we show that the radii of convexity of these functions are the smallest positive roots of some transcendental equations. Moreover, we will present some inequalities for the radii of convexity of the same functions.
\begin{theorem}
	Let ${\nu}>-1.$ Then the radius of convexity $r^{c}(g_\nu)$ of the function $$z\mapsto{g_\nu}(z)=2^{\nu}{\Gamma(\nu+1)}{z^{1-\nu}}{J_{\nu}(z)}$$ is the smallest positive root of the equation ${(z{g_{\nu}^{\prime}}(z))}^{\prime}=0$ and satisfies the following inequalities
	$$\frac{2\sqrt{\nu+1}}{3}<r^{c}(g_\nu)<6\sqrt{\frac{(\nu+1)(\nu+2)}{56\nu+137}},$$
	$$2\sqrt[4]{\frac{(\nu+1)^{2}(\nu+2)}{56\nu+137}}<r^{c}(g_\nu)<\sqrt{\frac{2(56\nu+137)(\nu+1)(\nu+3)}{208{\nu}^{2}+1172{\nu}+1693}},$$
	 $$\sqrt[6]{\frac{32(\nu+1)^{3}(\nu+2)(\nu+3)}{208{\nu}^{2}+1172{\nu}+1693}}<r^{c}(g_\nu)<2\sqrt{\frac{2(\nu+1)(\nu+2)(\nu+4)(208{\nu}^{2}+1172{\nu}+1693)}{3104{\nu}^{4}+36768{\nu}^{3}+161424{\nu}^{2}+312197{\nu}+223803}}.$$
\end{theorem}
\begin{theorem}
	Let ${\nu}>-1.$ Then the radius of convexity $r^{c}(h_\nu)$ of the function
	$$z\mapsto{h_\nu}(z)=2^{\nu}{\Gamma(\nu+1)}{z^{1-\frac{\nu}{2}}}{J_{\nu}(\sqrt{z})}$$ is the smallest positive root of the equation ${(z{h_{\nu}^{\prime}}(z))}^{\prime}=0$ and satisfies the following inequalities
	$${\nu}+1<r^{c}(h_\nu)<\frac{16(\nu+1)(\nu+2)}{7{\nu}+23},$$
	$$\sqrt{\frac{16(\nu+1)^{2}(\nu+2)}{7{\nu}+23}}<r^{c}(h_\nu)<\frac{2(\nu+1)(\nu+3)(7{\nu}+23)}{9{\nu}^{2}+60{\nu}+115},$$
	 $$\sqrt[3]{\frac{32(\nu+1)^{3}(\nu+2)(\nu+3)}{9{\nu}^{2}+60{\nu}+115}}<r^{c}(h_\nu)<\frac{8(\nu+1)(\nu+2)(\nu+4)(9{\nu}^{2}+60{\nu}+115)}{47{\nu}^{4}+621{\nu}^{3}+3136{\nu}^{2}+7221{\nu}+6195}.$$
\end{theorem}

\subsection{Radii of convexity of Struve functions}
In this subsection we consider two different normalized Struve functions of the first kind. Here we show that the radii of convexity of these functions are the smallest positive roots of some transcendental equations. We give also some lower and upper bounds for the radii of convexity of these functions.
\begin{theorem}
	Let $|\nu|\leq\frac{1}{2}$. Then the radius of convexity $r^{c}(u_\nu)$ of the function
	$$z\mapsto{u}_{\nu}(z)=\sqrt{\pi}2^{\nu}z^{-\nu}\Gamma\left({\nu}+\frac{3}{2}\right){\mathbf{H}}_{\nu}(z)$$ is the smallest positive root of the equation ${(z{u_{\nu}^{\prime}}(z))}^{\prime}=0$ and satisfies the following inequalities
	$$\sqrt{\frac{2{\nu}+3}{3}}<r^{c}(u_\nu)<\sqrt{\frac{36{\nu}^{2}+144{\nu}+135}{34{\nu}+105}},$$
	$$\sqrt[4]{\frac{3(2{\nu}+3)^{2}(2{\nu}+5)}{34{\nu}+105}}<r^{c}(u_\nu)<\sqrt{\frac{5(2{\nu}+3)(2{\nu}+7)(34{\nu}+105)}{3(268{\nu}^{2}+1824{\nu}+3213)}},$$
	 $$\sqrt[6]{\frac{5(2{\nu}+3)^{3}(2{\nu}+5)(2{\nu}+7)}{268{\nu}^{2}+1824{\nu}+3213}}<r^{c}(u_\nu)<3\sqrt{\frac{7(2{\nu}+3)(2{\nu}+5)(2{\nu}+9)(268{\nu}^{2}+1824{\nu}+3213)}{{\nu}^*}},$$ where ${\nu}^*=160336{\nu}^{4}+2256464{\nu}^{3}+11855904{\nu}^{2}+27626796{\nu}+24017715$.
\end{theorem}
\begin{theorem}
	Let $|\nu|\leq\frac{1}{2}$. Then the radius of convexity $r^{c}(w_\nu)$ of the function
	$$z\mapsto{w}_{\nu}(z)=\sqrt{\pi}2^{\nu}z^{\frac{1-\nu}{2}}\Gamma\left({\nu}+\frac{3}{2}\right){\mathbf{H}}_{\nu}(\sqrt{z})$$ is the smallest positive root of the equation ${(z{w_{\nu}^{\prime}}(z))}^{\prime}=0$ and satisfies the following inequalities
	$${\frac{3(2{\nu}+3)}{4}}<r^{c}(w_\nu)<{\frac{30(2{\nu}+3)(2{\nu}+5)}{26{\nu}+119}},$$
	$$\sqrt{\frac{45(2{\nu}+3)^{2}(2{\nu}+5)}{2(26{\nu}+119)}}<r^{c}(w_\nu)<{\frac{21(2{\nu}+3)(2{\nu}+7)(26{\nu}+119)}{2(404{\nu}^{2}+3396{\nu}+8665)}},$$
	 $$\sqrt[3]{\frac{945(2{\nu}+3)^{3}(2{\nu}+5)(2{\nu}+7)}{4(404{\nu}^{2}+3396{\nu}+8665)}}<r^{c}(w_\nu)<{\frac{30(2{\nu}+3)(2{\nu}+5)(2{\nu}+9)(404{\nu}^{2}+3396{\nu}+8665)}{{\nu}^{**}}},$$ where ${\nu}^{**}=36368{\nu}^{4}+588848{\nu}^{3}+3695776{\nu}^{2}+10793332{\nu}+11828151$.
\end{theorem}

\section{Proofs of the Main Results}
\setcounter{equation}{0}
\begin{proof}[Proof of Theorem 1]
	It is known (see \cite{bdoy}) that the zeros of the function $$h_{\nu}(z)=\sum_{n\geq0}\frac{(-1)^n{(2n+\nu+\alpha+1)}}{2^{2n}(\nu+\alpha+1)(\frac{3}{2})_{n}(\nu+\frac{3}{2})_{n}}z^{n}$$ all are real when $\alpha+\nu>-1$ and $\rvert\nu\rvert<\frac{1}{2}$. As a result of this we can say that the function $h_{\nu}$ belongs to the Laguerre-P\'{o}lya class $\mathcal{LP}$ of real entire functions, which are uniform limits of real polynomials whose all zeros are real. Thus, the function $z\mapsto h_{\nu}(z)$ has only real zeros and having growth order $\frac{1}{2}$ it can be written as the product
	$$h_{\nu}(z)=\prod_{n\geq 1}\left(1-\frac{z}{\zeta_{{\nu},n}^{2}}\right),$$
	where $\zeta_{{\nu},n}>0$ for each $n\in\mathbb{N}.$ By considering the Euler-Rayleigh sum $\delta_k=\sum_{n\geq1}{\zeta_{{\nu},n}^{-2k}}$ and the infinite sum representation of the Struve function $\mathbf{H}_{\nu}$ we have
	\begin{equation}
	 \frac{h^{\prime}_{\nu}(z)}{h_{\nu}(z)}=\sum_{n\geq1}\frac{1}{z-\zeta_{{\nu},n}^{2}}=-\sum_{k\geq0}\sum_{n\geq1}\frac{1}{(\zeta_{{\nu},n}^{2})^{k+1}}=-\sum_{k\geq0}\delta_{k+1}z^k, \rvert{z}\rvert<\zeta_{{\nu},1}^{2},
	\end{equation}
	\begin{equation}
	\frac{h^{\prime}_{\nu}(z)}{h_{\nu}(z)}=\sum_{n\geq0}\theta_nz^n\bigg/\sum_{n\geq0}\gamma_nz^n,
	\end{equation}
	where $$\theta_n=\frac{(-1)^{n+1}(2n+\nu+\alpha+3)(n+1)}{2^{2n+2}(\nu+\alpha+1)(\frac{3}{2})_{n+1}(\nu+\frac{3}{2})_{n+1}} \text{ and } \gamma_n=\frac{(-1)^{n}(2n+\nu+\alpha+1)}{2^{2n}(\nu+\alpha+1)(\frac{3}{2})_{n}(\nu+\frac{3}{2})_{n}}.$$ By comparing the coefficients of $(3.1)$ and $(3.2)$ we have the followings: $$\delta_1=\frac{(\alpha+\nu+3)}{3(2\nu+3)(\alpha+\nu+1)}, \delta_2=\frac{\kappa_1}{45(2\nu+3)^2(2\nu+5)(\alpha+\nu+1)^2},$$ $$\delta_3=\frac{\kappa_2}{945(2\nu+3)^3(4{\nu}^2+24\nu+35)(\alpha+\nu+1)^2},$$ where $$\kappa_1=-2\alpha^2\nu+7\alpha^2-4\alpha\nu^2+2\alpha\nu+42\alpha-2\nu^3-5\nu^2+72\nu+135$$ and
	\begin{align*}
	\kappa_2&=-4\alpha^3\nu^2-96\alpha^3\nu+145\alpha^3-12\alpha^2\nu^3-324\alpha^2\nu^2-429\alpha^2\nu+1305\alpha^2-12\alpha\nu^4-360\alpha\nu^3\\
	&-1689\alpha\nu^2+1170\alpha\nu+6291\alpha-4\nu^5-132\nu^4-1115\nu^3+621\nu^2+12339\nu+14931.
	\end{align*}
	Now by using the Euler-Rayleigh inequalities ${\delta_k}^{-\frac{1}{k}}<{{\zeta}_{\nu,1}}^2<\frac{\delta_k}{\delta_k+1}$ for $\alpha+\nu>-1, \rvert\nu\rvert<\frac{1}{2}$ and $k\in\{1,2,3\}$ we get the following lower bounds $$\zeta_{{\nu},1}^{2}>\frac{3(2{\nu}+3)(\alpha+{\nu}+1)}{\alpha+{\nu}+3},$$$$\zeta_{{\nu},1}^{2}>\frac{3(2{\nu}+3)(\alpha+{\nu}+1)\sqrt{5(2\nu+5)}}{\sqrt{\kappa_1}},$$$$\zeta_{{\nu},1}^{2}>\frac{3(2{\nu}+3)(\alpha+{\nu}+1)\sqrt[3]{35(2\nu+5)(2\nu+7)}}{\sqrt[3]{\kappa_2}}$$ and the upper bounds $$\zeta_{{\nu},1}^{2}<\frac{15(2\nu+3)(2\nu+5)(\alpha+\nu+1)(\alpha+\nu+3)}{\kappa_1},$$$$\zeta_{{\nu},1}^{2}<\frac{21(2\nu+3)(2\nu+7)(\alpha+\nu+1)\kappa_1}{\kappa_2}.$$
\end{proof}

\begin{proof}[Proof of Theorem 3]
	The normalized Lommel function $$l_{\mu}(z)=\frac{2\mu(\mu+1)}{(2\mu+1)}z^{-\frac{2\mu+1}{4}}\mathcal{L}_\mu(\sqrt{z})$$ has only real zeros for $\mu\in(-1,1), \mu\neq0$ and $\mu\neq-\frac{1}{2}$ (see \cite{BY}). Consequently, the function $l_{\mu}$ belongs to the Laguerre-P\'{o}lya class $\mathcal{LP}$ of real entire functions. Thus, $l_{\mu}(z)$ can be written as the product
	$$\prod_{n\geq1}\left(1-\frac{z}{{\tau_{\mu,n}^2}}\right)$$ where ${\tau}_{\mu,n}>0$ for each $n\in\mathbb{N}.$ Now by using the Euler-Rayleigh sum $\eta_k=\sum_{n\geq1}\tau_{\mu,n}^{-2k}$ and the infinite sum representation of the Lommel function $s_{\mu-\frac{1}{2},\frac{1}{2}}$ we get
	\begin{equation}
	 \frac{{l_\mu}^{\prime}(z)}{l_{\mu}(z)}=\sum_{n\geq1}\frac{1}{z-{\tau^2}_{\mu,n}}=-\sum_{n\geq1}\sum_{k\geq0}\frac{1}{{({\tau^2}_{\mu,n})^{k+1}}}z^k=-\sum_{k\geq0}\eta_{k+1}z^k, \rvert{z}\rvert<\tau_{\mu,1}^2,
	\end{equation}
	\begin{equation}
	\frac{{l_\mu}^{\prime}(z)}{l_{\mu}(z)}=\sum_{n\geq0}\rho_nz^n\bigg/\sum_{n\geq0}\sigma_nz^n,
	\end{equation}
	where $$\rho_n=\frac{(-1)^n(2n+\mu+\frac{5}{2})}{2^{2n+2}(\mu+\frac{1}{2})(\frac{\mu+2}{2})_{n+1}(\frac{\mu+3}{2})_{n+1}} \text{ and } \sigma_n=\frac{(-1)^n(2n+\mu+\frac{1}{2})}{2^{2n}(\mu+\frac{1}{2})(\frac{\mu+2}{2})_{n}(\frac{\mu+3}{2})_{n}}.$$ By equating the coefficients of $(3.3)$ and $(3.4)$ we obtain $$\eta_1=\frac{2\mu+5}{2\mu^3+11\mu^2+17\mu+6}, \eta_2=\frac{-4\mu^4-24\mu^3+19\mu^2+295\mu+392}{(\mu+2)^2(\mu+3)^2(\mu+4)(\mu+5)(2\mu+1)^2}$$ and
	$$\eta_3=\frac{8\mu^7+44\mu^6-554\mu^5-4731\mu^4-7672\mu^3+23551\mu^2+85834\mu+72384}{(\mu+2)^3(\mu+3)^3(\mu+4)(\mu+5)(\mu+6)(\mu+7)(2\mu+1)^3}.$$ Now by considering Euler-Rayleigh inequalities ${\eta_k}^{-\frac{1}{k}}<{\tau_{\mu,1}}^2<\frac{\eta_k}{\eta_{k+1}}$ for $\mu\in(-1,1), \mu\neq0, \mu\neq-\frac{1}{2}$ and $k\in\{1,2,3\}$ we obtain the following lower bounds 	$$(\tau_{\mu,1})^2>\frac{(\mu+2)(\mu+3)(2\mu+1)}{2\mu+5},$$
	$$(\tau_{\mu,1})^2>\frac{(\mu+2)(\mu+3)(2\mu+1)\sqrt{(\mu+4)(\mu+5)}}{\sqrt{-4\mu^4-24\mu^3+19\mu^2+295\mu+392}},$$
	 $$(\tau_{\mu,1})^2>\frac{(\mu+2)(\mu+3)(2\mu+1)\sqrt[3]{(\mu+4)(\mu+5)(\mu+6)(\mu+7)}}{\sqrt[3]{8\mu^7+44\mu^6-554\mu^5-4731\mu^4-7672\mu^3+23551\mu^2+85834\mu+72384}}$$ and the upper bounds
	$$(\tau_{\mu,1})^2<\frac{(\mu+2)(\mu+3)(\mu+4)(\mu+5)(2\mu+1)(2\mu+5)}{-4\mu^4-24\mu^3+19\mu^2+295\mu+392},$$
	 $$(\tau_{\mu,1})^2<\frac{(\mu+2)(\mu+3)(\mu+6)(\mu+7)(2\mu+1)(-4\mu^4-24\mu^3+19\mu^2+295\mu+392)}{8\mu^7+44\mu^6-554\mu^5-4731\mu^4-7672\mu^3+23551\mu^2+85834\mu+72384}.$$
\end{proof}
\begin{proof}[Proof of Theorem 4]
	If we consider the Maclaurin series expansion of the function $$z\mapsto{u}_{\nu}(z)=\left(\sqrt{\pi}2^{\nu}\Gamma\left(\nu+\frac{3}{2}\right)\mathbf{H}_{\nu}(z)\right)^{\frac{1}{\nu+1}}$$ we obtain
	\begin{equation}
	u(z)=z-\frac{1}{3(\nu+1)(2\nu+3)}z^3+\frac{1}{90(\nu+1)^2(2\nu+3)^2(2\nu+5)}z^5-{\dots}.
	\end{equation}
	Therefore, the function $u_{\nu}$ has real coefficients. Also we know that if the function $z\mapsto{z+\alpha_{2}z+{\dots}}$ has real coefficients, then its radius of starlikeness are less or equal than its radius of univalence, see \cite{wilf}. Now, we should show that the radii of univalence are less or equal than the corresponding radii of starlikeness. From the definition of $u_{\nu}(z)$ we can write that
	\begin{equation}
	 \frac{zu^{\prime}_{\nu}(z)}{u_{\nu}(z)}=\frac{1}{\nu+1}\frac{z\mathbf{H}^{\prime}_{\nu}(z)}{\mathbf{H}_{\nu}(z)}=1-\frac{2}{\nu+1}\sum_{n\geq1}\frac{z^2}{h^{2}_{\nu,n}-z^2}.
	\end{equation}
	Thus, for $\nu\in[-\frac{1}{2},\frac{1}{2}],$ we obtain that
	\begin{align*}
	 \real\left(\frac{zu^{\prime}_{\nu}(z)}{u_{\nu}(z)}\right)=1-\frac{2}{\nu+1}\sum_{n\geq1}\real\left(\frac{z^2}{h^{2}_{\nu,n}-z^2}\right)\geq1-\frac{2}{\nu+1}\sum_{n\geq1}\frac{|z|^2}{h^{2}_{\nu,n}-|z|^2}=\frac{|z|u^{\prime}_{\nu}(|z|)}{u_{\nu}(|z|)}.
	\end{align*}
	That is,
	\begin{equation}
	\real\left(\frac{zu^{\prime}_{\nu}(z)}{u_{\nu}(z)}\right)>\frac{ru^{\prime}_{\nu}(r)}{u_{\nu}(r)},
	\end{equation}
	where $r=\left|z\right|.$
	The quantity on the right-hand side of the inequality $(3.7)$ remains positive until the first positive zero of $u^{\prime}_{\nu}.$ These show that indeed the radius of univalence corresponds to the radius of starlikeness of the function $u_{\nu}$.
	
\end{proof}
\begin{proof}[Proof of Theorem 5]
	If we consider the Maclaurin series expansion of the function
	$$z\mapsto{f_{\mu}}(z)={f_{\mu-\frac{1}{2},\frac{1}{2}}}(z)=\left(\mu(\mu+1)s_{\mu-\frac{1}{2},\frac{1}{2}}(z)\right)^\frac{1}{\mu+\frac{1}{2}}$$
	we obtain
	$${f_{\mu}}(z)=z-\frac{2}{(\mu+2)(\mu+3)(2\mu+1)}z^3+\frac{2{\mu}^3+16{\mu}^2+39\mu-16}{2(\mu+2)(\mu+3)(\mu+4)(\mu+5)(2\mu+1)^2}z^5-{\dots}.$$
	Therefore the radius of starlikeness of the function $f_{\mu}$ is less or equal than its radius of univalence, see \cite{wilf}. On the other hand, from the definition of $f_{\mu}$ we can write that
	\begin{equation}
	 \frac{z{f^{\prime}_{\mu}}(z)}{{f_{\mu}}(z)}=\frac{1}{1+\frac{\mu}{2}}\frac{zs^{\prime}_{\mu-\frac{1}{2},\frac{1}{2}}(z)}{s_{\mu-\frac{1}{2},\frac{1}{2}}(z)}=1-\frac{2}{1+\frac{\mu}{2}}\sum_{n\geq1}\frac{z^2}{l^{2}_{\mu,n}-z^2}.
	\end{equation}
	Thus, for $\mu\in(-\frac{1}{2},1), \mu\neq0$ we obtain that
	\begin{align*}
	 \real\left(\frac{zf^{\prime}_{\mu}(z)}{f_{\mu}(z)}\right)=1-\frac{2}{1+\frac{\mu}{2}}\sum_{n\geq1}\real\left(\frac{z^2}{l^{2}_{\mu,n}-z^2}\right)\geq1-\frac{2}{{1+\frac{\mu}{2}}}\sum_{n\geq1}\frac{|z|^2}{l^{2}_{\mu,n}-|z|^2}=\frac{|z|f^{\prime}_{\mu}(|z|)}{f_{\mu}(|z|)}.
	\end{align*}
	That is
	\begin{equation}
	\real\left(\frac{zf^{\prime}_{\mu}(z)}{f_{\mu}(z)}\right)>\frac{rf^{\prime}_{\mu}(r)}{f_{\mu}(r)}
	\end{equation}
	where $r=\left|z\right|.$
	The quantity on the right-hand side of the inequality $(3.9)$ remains positive until the first positive zero of $f^{\prime}_{\mu}$ is reached. These show that indeed the radius of univalence corresponds to the radius of starlikeness of the function $f_{\mu}$.
\end{proof}

\begin{proof}[Proof of Theorem 6]
	By using the Alexander duality theorem for starlike and convex functions we can say that the function $g_{\nu}$ is convex if and only if $z\mapsto{zg^{\prime}_{\nu}(z)}$ is starlike. But, the smallest positive zero of $z\mapsto{\left(zg^{\prime}_{\nu}(z)\right)^{\prime}}$ is actually the radius of starlikeness of $z\mapsto{zg^{\prime}_{\nu}(z)}$, according to \cite{bsk,bos}. Therefore, the radius of convexity $r^{c}(g_\nu)$ is the smallest positive root of the equation ${(z{g_{\nu}^{\prime}}(z))}^{\prime}=0$. See also \cite{basz} for more details. Now, by considering the Bessel differential equation
	\begin{equation}
	z^2J^{\prime \prime}_{\nu}(z)+zJ^{\prime}_{\nu}(z)+(z^{2}-{\nu}^2)J_{\nu}(z)=0
	\end{equation}
	 and the infinite series representations of Bessel function and its derivative
	 \begin{equation}
	 J_{\nu}(z)=\sum_{n\geq0}\frac{(-1)^n{z^{2n+\nu}}}{2^{2n+\nu}n!\Gamma(n+\nu+1)},
	 \end{equation}
	 \begin{equation}
	 J^{\prime}_{\nu}(z)=\sum_{n\geq0}\frac{(-1)^n(2n+\nu){z^{2n+\nu-1}}}{2^{2n+\nu}n!\Gamma(n+\nu+1)},
	 \end{equation}
	respectively, we obtain
	\begin{equation}
	\Delta_{\nu}(z)=\left(zg^{\prime}_{\nu}(z)\right)^{\prime}=1+\sum_{n\geq1}\frac{(-1)^n(2n+1)^2{z^{2n}}}{2^{2n}{n!}{(\nu+1)_n}}.
	\end{equation}
	 Since the function $g_{\nu}$ belongs to the Laguerre-P\'olya class of entire functions and $\mathcal{LP}$ is closed under differentiation, we can say that the function $\Delta_{\nu}$ belongs also to the Laguerre-P\'olya class. Therefore, the zeros of the function $\Delta_{\nu}$ are all real. Suppose that $\beta_{\nu,n}$'s are the zeros of the function $\Delta_{\nu}$. Then the function $\Delta_{\nu}$ has the infinite product representation as follows:
	\begin{equation}
	\Delta_{\nu}(z)=\prod_{n\geq1}\left(1-\frac{z^2}{\beta_{\nu,n}^2}\right).
	\end{equation}
	By taking the logarithmic derivative of $(3.14)$ we get
	\begin{equation}
	\frac{\Delta^{\prime}_{\nu}(z)}{\Delta_{\nu}(z)}=-2\sum_{k\geq0}\rho_{k+1}z^{2k+1}, \rvert{z}\rvert<\beta_{\nu,1}^2,
	\end{equation}
	where $\rho_{k}=\sum_{n\geq1}{\beta_{\nu,n}^{-2k}}$. On the other hand, by considering infinite sum representation of $\Delta_{\nu}(z)$ we obtain
	\begin{equation}
	\frac{\Delta^{'}_{\nu}(z)}{\Delta_{\nu}(z)}=\sum_{n\geq0}\xi_{n}z^{2n+1}\bigg/\sum_{n\geq0}\kappa_{n}z^{2n},
	\end{equation}
	where $$\xi_{n}=\frac{(-1)^{n+1}2(2n+3)^{2}}{2^{2n+2}{n!}(\nu+1)_{n+1}} \text{ and } \kappa_{n}=\frac{(-1)^{n}(2n+1)^{2}}{2^{2n}{n!}(\nu+1)_{n}}.$$ By comparing the coefficients of $(3.15)$ and $(3.16)$ we obtain $$\rho_{1}=\frac{9}{4(\nu+1)}, \rho_{2}=\frac{56\nu+137}{16(\nu+1)^2(\nu+2)}, \rho_{3}=\frac{208{\nu}^2+1172{\nu}+1693}{32(\nu+1)^3(\nu+2)(\nu+3)}$$ and $$\rho_{4}=\frac{3104{\nu}^4+36768{\nu}^3+161424{\nu}^2+312197{\nu}+223803}{216(\nu+1)^4(\nu+2)^2(\nu+3)(\nu+4)}.$$
	Now by considering the Euler-Rayleigh inequalities ${\rho_k}^{-\frac{1}{k}}<{{\beta}_{\nu,1}}^2<\frac{\rho_k}{\rho_{k+1}}$ for $\nu>-1$ and $k\in\{1,2,3\}$ we obtain following inequalities
	$$\frac{2\sqrt{\nu+1}}{3}<r^{c}(g_\nu)<6\sqrt{\frac{(\nu+1)(\nu+2)}{56\nu+137}},$$
	$$2\sqrt[4]{\frac{(\nu+1)^{2}(\nu+2)}{56\nu+137}}<r^{c}(g_\nu)<\sqrt{\frac{2(56\nu+137)(\nu+1)(\nu+3)}{208{\nu}^{2}+1172{\nu}+1693}},$$
	 $$\sqrt[6]{\frac{32(\nu+1)^{3}(\nu+2)(\nu+3)}{208{\nu}^{2}+1172{\nu}+1693}}<r^{c}(g_\nu)<2\sqrt{\frac{2(\nu+1)(\nu+2)(\nu+4)(208{\nu}^{2}+1172{\nu}+1693)}{3104{\nu}^{4}+36768{\nu}^{3}+161424{\nu}^{2}+312197{\nu}+223803}}.$$

\end{proof}
\begin{proof}[Proof of Theorem 7]
	By using the same procedure as in the previous proof we can say that the radius of convexity $r^{c}(h_\nu)$ is the smallest positive root of the equation ${(z{h_{\nu}^{\prime}}(z))}^{\prime}=0$. See also \cite{basz} for more details. Now, by setting $\sqrt{z}$ instead of $z$ in the $(3.10), (3.11) \text{ and }(3.12),$ respectively, we obtain
	\begin{equation}
	\theta_{\nu}(z)={(z{h_{\nu}^{\prime}}(z))}^{\prime}=1+\sum_{n\geq1}\frac{(-1)^n(n+1)^2{z^{n}}}{2^{2n}{n!}{(\nu+1)_n}}.
	\end{equation}
	In addition, we know that $h_{\nu}$ belongs to the Laguerre-P\'olya class of entire functions $\mathcal{LP}$. Since $\mathcal{LP}$ is closed under differentiation, we can say that the function $\theta_{\nu}$ belongs also to the Laguerre-P\'olya class. That is, the zeros of the function $\theta_{\nu}$ are all real. Suppose that $\gamma_{\nu,n}$'s are the zeros of the function $\theta_{\nu}$. Then the function $\theta_{\nu}$ has the infinite product representation as follows:
	\begin{equation}
	\theta_{\nu}(z)=\prod_{n\geq1}\left(1-\frac{z}{\gamma_{\nu,n}}\right).
	\end{equation}
	By logarithmic derivation of $(3.18)$ we get
	\begin{equation}
	\frac{\theta^{\prime}_{\nu}(z)}{\theta_{\nu}(z)}=-\sum_{k\geq0}\varrho_{k+1}z^{k},
	\end{equation}
	where $\varrho_{k}=\sum_{n\geq1}{\gamma_{\nu,n}^{-k}}.$ Also, by using derivative of infinite sum representation of $\theta_{\nu}(z)$ we get
	\begin{equation}
	\frac{\theta^{\prime}_{\nu}(z)}{\theta_{\nu}(z)}=\sum_{n\geq0}m_{n}z^{n}\bigg/\sum_{n\geq0}s_{n}z^{n}, \left|z\right|<\gamma_{\nu,1},
	\end{equation}
	where $$m_{n}=\frac{(-1)^{n+1}(n+2)^2}{2^{2n}{n!}(\nu+1)_{n+1}} \text{ and } s_{n}=\frac{(-1)^{n}(n+1)^2}{2^{2n}{n!}(\nu+1)_{n}}.$$ By comparing the coefficients of $(3.19)$ and $(3.20)$ we have $$\varrho_{1}=\frac{1}{\nu+1}, \varrho_{2}=\frac{7\nu+23}{16(\nu+1)^2(\nu+2)}, \varrho_{3}=\frac{9\nu^2+60\nu+115}{32(\nu+1)^3(\nu+2)(\nu+3)}$$ and $$\varrho_{4}=\frac{47\nu^4+621\nu^3+3136\nu^2+7221\nu+6195}{256(\nu+1)^4(\nu+2)^2(\nu+3)(\nu+4)}.$$ By applying the Euler-Rayleigh inequalities ${\varrho_k}^{-\frac{1}{k}}<{{\gamma}_{\nu,1}}<\frac{\varrho_k}{\varrho_{k+1}}$ for $\nu>-1$ and $k\in\{1,2,3\}$ we have
	$${\nu}+1<r^{c}(h_\nu)<\frac{16(\nu+1)(\nu+2)}{7{\nu}+23},$$
	$$\sqrt{\frac{16(\nu+1)^{2}(\nu+2)}{7{\nu}+23}}<r^{c}(h_\nu)<\frac{2(\nu+1)(\nu+3)(7{\nu}+23)}{9{\nu}^{2}+60{\nu}+115},$$
	 $$\sqrt[3]{\frac{32(\nu+1)^{3}(\nu+2)(\nu+3)}{9{\nu}^{2}+60{\nu}+115}}<r^{c}(h_\nu)<\frac{8(\nu+1)(\nu+2)(\nu+4)(9{\nu}^{2}+60{\nu}+115)}{47{\nu}^{4}+621{\nu}^{3}+3136{\nu}^{2}+7221{\nu}+6195}.$$
\end{proof}
\begin{proof}[Proof of Theorem 8]
	Similarly as in the proof of Theorem 6 we observe that the radius of convexity $r^{c}(u_\nu)$ is the smallest positive root of the equation ${(z{u_{\nu}^{\prime}}(z))}^{\prime}=0$. See also \cite{BY} fore more details. Now, by considering the Struve differential equation
	\begin{equation}
	z^2\mathbf{H}_{\nu}^{\prime \prime}(z)+z\mathbf{H}_{\nu}^{\prime}(z)+(z^2-\nu^2)\mathbf{H}_{\nu}(z)=\frac{4(\frac{z}{2})^{\nu+1}}{\sqrt{\pi}\Gamma\left(\nu+\frac{1}{2}\right)}
	\end{equation}
	and the infinite series representations of Struve function and its derivative
	\begin{equation}
	\mathbf{H}_{\nu}(z)=\sum_{n\geq0}\frac{(-1)^n}{\Gamma\left(n+\frac{3}{2}\right)\Gamma\left(\nu+n+\frac{3}{2}\right)}\left(\frac{z}{2}\right)^{2n+\nu+1},
	\end{equation}
	\begin{equation}
	 \mathbf{H}_{\nu}^{\prime}(z)=\sum_{n\geq0}\frac{(-1)^n(2n+\nu+1)}{2\Gamma\left(n+\frac{3}{2}\right)\Gamma\left(\nu+n+\frac{3}{2}\right)}\left(\frac{z}{2}\right)^{2n+\nu},
	\end{equation}
	respectively, we get
	\begin{equation}
	\Omega_{\nu}(z)={(z{u_{\nu}^{\prime}}(z))}^{\prime}=1+\sum_{n\geq1}\frac{(-1)^n(2n+1)}{2^{2n}(\frac{1}{2})_{n}(\nu+\frac{3}{2})_{n}}z^{2n}.
	\end{equation}
	Since the function ${u_{\nu}}$ belongs to the Laguerre-P\'olya class of entire functions $\mathcal{LP}$ and this class is closed under differentation we obtain that the function $\Omega_{\nu}$ belongs also to the Laguerre-P\'olya class. Therefore, the zeros of the function $\Omega_{\nu}$ are all real. Suppose that $\vartheta_{\nu,n}$'s are the zeros of the function $\Omega_{\nu}$. Then the function $\Omega_{\nu}$ has infinite product representation as follows:
	\begin{equation}
	\Omega_{\nu}(z)=\prod_{n\geq1}\left(1-\frac{z^2}{\vartheta_{\nu,n}^2}\right).
	\end{equation}
	By taking the logarithmic derivative of $(3.25)$ we have
	\begin{equation}
	\frac{\Omega^{\prime}_{\nu}(z)}{\Omega_{\nu}(z)}=-2\sum_{k\geq0}\chi_{k+1}z^{2k+1}, \left|z\right|<{\vartheta^{2}_{\nu,1}},
	\end{equation}
	where $\chi_{k}=\sum_{n\geq1}{\vartheta_{\nu,n}^{-2k}}.$ On the other hand, by considering infinite sum representation of $\Omega_{\nu}(z)$ we get
	\begin{equation}
	\frac{\Omega^{\prime}_{\nu}(z)}{\Omega_{\nu}(z)}=\sum_{n\geq0}\tau_{n}z^{2n+1}\bigg/\sum_{n\geq0}\varsigma_{n}z^{2n},
	\end{equation}
	where $$\tau_n=\frac{(-1)^{n+1}(2n+3)(n+1)}{2^{2n+1}(\frac{1}{2})_{n+1}(\nu+\frac{3}{2})_{n+1}} \text{ and } \varsigma_n=\frac{(-1)^{n}(2n+1)}{2^{2n}(\frac{1}{2})_{n}(\nu+\frac{3}{2})_{n}}.$$ Now, by comparing the coefficients of $(3.26) \text{ and }(3.27)$ we obtain $$\chi_1=\frac{3}{2\nu+3}, \chi_2=\frac{34\nu+105}{3(2\nu+3)^2(2\nu+5)}, \chi_3=\frac{268\nu^2+1824\nu+3213}{5(2\nu+3)^3(2\nu+5)(2\nu+7)}$$ and $$\chi_4=\frac{160336\nu^4+2256464\nu^3+11855904\nu^2+27626796\nu+24017715}{315(2\nu+3)^4(2\nu+5)^2(2\nu+7)(2\nu+9)}.$$
	By using the Euler-Rayleigh inequalities ${\chi_k}^{-\frac{1}{k}}<{{\vartheta}_{\nu,1}}^2<\frac{\chi_k}{\chi_k+1}$ for $|\nu|\leq\frac{1}{2}$  and $k\in\{1,2,3\}$ we obtain
	$$\sqrt{\frac{2{\nu}+3}{3}}<r^{c}(u_\nu)<\sqrt{\frac{36{\nu}^{2}+144{\nu}+135}{34{\nu}+105}},$$
	$$\sqrt[4]{\frac{3(2{\nu}+3)^{2}(2{\nu}+5)}{34{\nu}+105}}<r^{c}(u_\nu)<\sqrt{\frac{5(2{\nu}+3)(2{\nu}+7)(34{\nu}+105)}{3(268{\nu}^{2}+1824{\nu}+3213)}},$$
	 $$\sqrt[6]{\frac{5(2{\nu}+3)^{3}(2{\nu}+5)(2{\nu}+7)}{268{\nu}^{2}+1824{\nu}+3213}}<r^{c}(u_\nu)<3\sqrt{\frac{7(2{\nu}+3)(2{\nu}+5)(2{\nu}+9)(268{\nu}^{2}+1824{\nu}+3213)}{{\nu}^*}},$$ where ${\nu}^*=160336{\nu}^{4}+2256464{\nu}^{3}+11855904{\nu}^{2}+27626796{\nu}+24017715$.
\end{proof}
\begin{proof}[Proof of Theorem 9]
	By using the same idea as in the proof of Theorem 6 we have that  the radius of convexity $r^{c}(w_\nu)$ is the smallest positive root of the equation ${(z{w_{\nu}^{\prime}}(z))}^{\prime}=0$. See also \cite{BY} fore more details. Now, if we put $\sqrt{z}$ instead of $z$ in the $(3.21), (3.22) \text{ and }(3.23)$ respectively, after some calculations we obtain
	\begin{equation}
	\psi_{\nu}(z)={(z{w_{\nu}^{\prime}}(z))}^{\prime}=1+\sum_{n\geq1}\frac{(-1)^n(n+1)^2}{2^{2n}(2n+1)(\frac{1}{2})_{n}(\nu+\frac{3}{2})_{n}}z^{n}.
	\end{equation}
	On the other hand, we know that the function ${w_{\nu}}$ belongs to the Laguerre-P\'olya class of entire functions $\mathcal{LP}$ and the Laguerre-P\'olya class of entire functions is closed under differentiation. Therefore, we get that the function $\psi_{\nu}$ belongs also to the Laguerre-P\'olya class. Hence, the zeros of the function $\psi_{\nu}$ are all real. Suppose that $\epsilon_{\nu,n}$'s are the zeros of the function $\psi_{\nu}$. Then the function $\psi_{\nu}$ has the infinite product representation as follows:
	\begin{equation}
	\psi_{\nu}(z)=\prod_{n\geq1}\left(1-\frac{z}{\epsilon_{\nu,n}}\right).
	\end{equation}
	If we take the derivative of the $(3.29)$ logarithmically then we get
	\begin{equation}
	\frac{\psi^{\prime}_{\nu}(z)}{\psi_{\nu}(z)}=-\sum_{k\geq0}\varphi_{k+1}z^{k}, \left|z\right|<{\epsilon_{\nu,1}},
	\end{equation}
	where $\varphi_{k}=\sum_{n\geq1}{\epsilon_{\nu,n}^{-k}}$. Also, by taking derivative of $(3.28)$ we have
	\begin{equation}
	\frac{\psi^{\prime}_{\nu}(z)}{\psi_{\nu}(z)}=\sum_{n\geq0}t_{n}z^{n}\bigg/\sum_{n\geq0}r_{n}z^{n},
	\end{equation}
	where $$t_n=\frac{(-1)^{n+1}(n+2)^2(n+1)}{2^{2n+2}(2n+3)(\frac{1}{2})_{n+1}(\nu+\frac{3}{2})_{n+1}} \text{ and } r_n=\frac{(-1)^{n}(n+1)^2}{2^{2n}(2n+1)(\frac{1}{2})_{n}(\nu+\frac{3}{2})_{n}}.$$
	Now, by comparing the coefficients of $(3.30)$ and $(3.31)$ we get $$\varphi_{1}=\frac{4}{3(2\nu+3)}, \varphi_{2}=\frac{2(26\nu+119)}{45(2\nu+3)^2(2\nu+5)}, \varphi_{3}=\frac{4(404\nu^2+3396\nu+8665)}{945(2\nu+3)^3(2\nu+5)(2\nu+7)}$$ and $$\varphi_{4}=\frac{2(36368\nu^4+588848\nu^3+3695776\nu^2+10793332\nu+11828151)}{14175(2\nu+3)^4(2\nu+5)^2(2\nu+7)(2\nu+9)}.$$ When we use the Euler-Rayleigh inequalities ${\varphi_k}^{-\frac{1}{k}}<{{\epsilon}_{\nu,1}}<\frac{\varphi_k}{\varphi_{k+1}}$ for  $|\nu|\leq\frac{1}{2}$ and $k\in\{1,2,3\}$ we obtain the following inequalities
	$${\frac{3(2{\nu}+3)}{4}}<r^{c}(w_\nu)<{\frac{30(2{\nu}+3)(2{\nu}+5)}{26{\nu}+119}},$$
	$$\sqrt{\frac{45(2{\nu}+3)^{2}(2{\nu}+5)}{2(26{\nu}+119)}}<r^{c}(w_\nu)<{\frac{21(2{\nu}+3)(2{\nu}+7)(26{\nu}+119)}{2(404{\nu}^{2}+3396{\nu}+8665)}},$$
	 $$\sqrt[3]{\frac{945(2{\nu}+3)^{3}(2{\nu}+5)(2{\nu}+7)}{4(404{\nu}^{2}+3396{\nu}+8665)}}<r^{c}(w_\nu)<{\frac{30(2{\nu}+3)(2{\nu}+5)(2{\nu}+9)(404{\nu}^{2}+3396{\nu}+8665)}{{\nu}^{**}}},$$ where ${\nu}^{**}=36368{\nu}^{4}+588848{\nu}^{3}+3695776{\nu}^{2}+10793332{\nu}+11828151$.
\end{proof}

\end{document}